\documentclass[12pt]{article}
\usepackage{amsthm,amsfonts,amssymb,amscd}

\begin{document}

\begin{center}
\large \bf Birational geometry of Fano double covers
\end{center}\vspace{0.5cm}

\centerline{A.V.Pukhlikov}\vspace{0.5cm}

\parshape=1
3cm 10cm \noindent {\small \quad\quad\quad \quad\quad\quad\quad
\quad\quad\quad {\bf }\newline We prove divisorial canonicity of
Fano double hypersurfaces of general position.} \vspace{1.5cm}

\section*{Introduction}\vspace{0.1cm}

{\bf 0.1. The main result.} The symbol ${\mathbb P}$ stands for
the projective space ${\mathbb P}^{M+1}$, $M\geq 6$. Fix a pair of
integers $m\geq 3$, $l\geq 2$, satisfying the relation $m+l=M+1$.
Let
$$
{\cal F}\subset\{(f,g)\in H^0({\mathbb P},{\cal O}_{\mathbb
P}(m))\times H^0({\mathbb P},{\cal O}_{\mathbb P}(2l))\}
$$
be the Zariski open set of pairs of non-zero polynomials $(f,g)$,
for which the double cover $\sigma\colon V\to Q\subset{\mathbb
P}$, branched over $W=W^*\cap Q$, is an irreducible variety, where
$Q=\{f=0\}$ is a hypersurface of dergee $m$, $W^*=\{g=0\}$ is a
hypersurface of degree $2l$. Set ${\cal F}_{\rm sm}\subset{\cal
F}$ to be the open subset, corresponding to smooth double covers
$V$. The aim of this paper is to prove the following fact.

{\bf Theorem 1.} {\it There exists a non-empty Zariski open subset
${\cal F}_{\rm reg}\subset {\cal F}_{\rm sm}$ such that for any
pair $(f,g)\in{\cal F}_{\rm reg}$ the corresponding variety $V$
satisfies the property of divisorial canonicity: for any effective
divisor $D\in|-nK_V|$, $n\geq 1$, the pair
$$
(V,\frac{1}{n}D)
$$
has canonical singularities.}

Recall that for $(f,g)\in{\cal F}_{\rm sm}$ the corresponding
variety $V=V(f,g)$ is a primitive Fano variety of dimension $M\geq
6$, that is, ${\mathop{\rm Pic}}V={\mathbb Z}K_V$. Canonicity of
the pair $(V,\frac{1}{n}D)$ means that for any birational morphism
$\varphi\colon V^+\to V$ and any exceptional divisor $E\subset
V^+$ the inequality
$$
\nu_E(D)\leq na(E)
$$
holds, where $a(E)$ is the discrepancy of $E$ with respect to the
model $V$, that is, the inequality, opposite to the Noether-Fano
inequality. By linearity of this inequality one may always assume
that $D$ is a prime divisor, that is, irreducible and
reduced.\vspace{0.3cm}

{\bf 0.2. Birational rigidity.} The property of divisorial
canonicity (the property ($C$)) was introduced in [1]. If the pair
$(V,\frac{1}{n}D)$ is canonical for a general divisor
$D\in\Sigma\subset|-nK_V|$ of any movable linear system $\Sigma$,
then the variety $V$ satisfies the property of {\it movable}
canonicity (the property ($M$)). Finally, if the pair
$(V,\frac{1}{n}D)$ is log canonical for any divisor
$D\in\Sigma\subset|-nK_V|$, then the variety $V$ is {\it
divisorially log canonical} (satisfies the property ($L$)). The
movable canonicity is shown for many classes of Fano varieties,
see [2] and the bibliography for that paper. The latter property
is important because it immediately implies {\it birational
rigidity} of the given variety.

Recall that a smooth projective rationally connected variety $X$
is said to be {\it birationally superrigid}, if for any movable
linear system $\Sigma$ on $X$ the equality
$$
c_{\rm virt}(\Sigma)=c(\Sigma,X)
$$
holds, where $c(\Sigma,X)=\mathop{\rm sup}\{t\in{\mathbb
Q}|D+tK_X\in A^1_+X,D\in\Sigma\}$ is the {\it threshold of
canonical adjunction} (the symbol $A^1_+X$ stands for the
pseudoeffective cone of the variety $X$ in $A^1_{\mathbb
R}X=A^1X\otimes{\mathbb R}$), whereas
$$
c_{\rm virt}(\Sigma)=\mathop{\rm inf}\limits_{X^+\to
X}\{c(\Sigma^+,X^+)\}
$$
is the {\it virtual threshold of canonical adjunction}, the
infimum is taken over all birational morphisms $X^+\to X$ of
smooth projective varieties, $\Sigma^+$ is the strict transform of
$\Sigma$ on $X^+$. If for any movable system $\Sigma$ there is a
birational self-map $\chi\in\mathop{\rm Bir}X$, providing the
equality of the thresholds,
$$
c_{\rm virt}(\Sigma)=c(\chi^{-1}_*\Sigma,X),
$$
where $\chi^{-1}_*\Sigma$ is the strict transform of the system
$\Sigma$ with respect to $\chi$, then the variety $X$ is said to
be {\it birationally rigid}. The property ($M$) and birational
superrigidity of Fano double hypersurfaces were proven in [3].

The main geometric implication of birational rigidity for
primitive Fano varieties is the absence of non-trivial structures
of a rationally connected fiber space, that is, of rational
dominant maps $\rho\,\colon X\dashrightarrow S$, the fiber of
general position of which is rationally connected. For this
reason, birationally rigid primitive Fano varieties are
automatically non-rational. Besides, birational rigidity makes it
possible to give an exhaustive description of birational maps of
the given variety onto other varieties.\vspace{0.3cm}

{\bf 0.3. The theorem on direct products.} The importance of
divisorial canonicity is connected with the following theorem
proven in [1] (the condition ($C$) implies the conditions ($M$)
and ($L$) in an obvious way).

{\bf Theorem 2.} {\it Assume that primitive Fano varieties
$F_1,\dots,F_K$, $K\geq 2$, satisfy the conditions ($L$) and
($M$). Then their direct product
$$
V=F_1\times\dots\times F_K
$$
is a birationally superrigid variety.}

Here are the main geometric consequences of birational rigidity of
the direct product $V$ (see [1, Corollary 1]).

(i) All structures of a rationally connected fiber space on the
variety $V$ are projections onto direct factors. More precisely,
let $\beta\colon V^{\sharp}\to S^{\sharp}$ be a rationally
connected fiber space and $\chi\colon V-\,-\,\to V^{\sharp}$ a
birational map. Then there exist a set of indices
$$
I=\{i_1,\dots,i_k\}\subset \{1,\dots,K\}
$$
and a birational map
$$
\alpha\colon F_I=\prod\limits_{i\in I}F_i-\,-\,\to S^{\sharp}
$$
such that the following diagram commutes:
$$
\begin{array}{rcccl}
& V &\stackrel{\chi}{-\,-\,\to} & V^{\sharp}&\\
\pi_I &\downarrow & &\downarrow &\beta\\
& F_I & \stackrel{\alpha}{-\,-\,\to}& S^{\sharp},&
\end{array}
$$
that is, $\beta\circ\chi=\alpha\circ\pi_I$, where
$\pi_I\colon\prod\limits^K_{i=1}F_i\to \prod\limits_{i\in I}F_i$
is the natural projection onto a direct factor. In particular, on
the variety $V$ there are no structures of a fibration into
rationally connected varieties of dimension strictly smaller than
$\mathop{\rm min}\{\mathop{\rm dim}F_i\}$. In particular, $V$ has
no structures of a fibration into conics and rational surfaces.

(ii) The groups of birational and biregular self-maps of the
variety $V$ coincide:
$$
\mathop{\rm Bir}V=\mathop{\rm Aut}V.
$$
In particular, the group $\mathop{\rm Bir}V$ is finite.

(iii) The variety $V$ is non-rational.

For a generic double cover $F\in{\cal F}_{\rm reg}$ we have (see
[3]):
$$
\mathop{\rm Bir}F={\mathop{\rm Aut}}F={\mathbb Z}/2{\mathbb Z}.
$$
Therefore, for pair-wise non-isomorphic double covers of general
position $F_1,\dots,F_K\in{\cal F}_{\rm reg}$ we get
$$
\mathop{\rm Bir}V=\mathop{\rm Aut}V=({\mathbb Z}/2{\mathbb Z})^K.
$$
On the other hand, for a double cover $F\in{\cal F}_{\rm reg}$ of
general position the group $\mathop{\rm Bir}F^{\times
K}={\mathop{\rm Aut}}F^{\times K}$ is the extension
$$
1\to({\mathbb Z/2{\mathbb Z}})^K\to\mathop{\rm Aut}F^{\times K}\to
S_K\to 1,
$$
where $S_K$ is the symmetric group of permutations of $K$
elements. This extension is in fact a semi-direct product: the
group of automorphisms $\mathop{\rm Aut}V$ contains $S_K$ as a
subgroup permuting the factors of the direct product
$F\times\dots\times F$. The action of $S_K$ on the subgroup
$({\mathbb Z}/2{\mathbb Z})^K$ is also obvious: the permutations
$\mu\in S_K$ permute the generators $\sigma_1,\dots,\sigma_K$ of
the Galois groups of the factors $F$.

Another application of Theorem 1 is in that it allows to apply the
linear method of proving birational rigidity [4] to investigating
Fano fibrations $V/{\mathbb P}^1$ into double hypersurfaces [5,6],
see \S 5 of this paper.\vspace{0.3cm}

{\bf 0.4. The structure of the paper.} The open set ${\cal F}_{\rm
reg}$ is defined by explicit conditions in \S 1. In \S 4 we show
that this set is non-empty. Let $(f,g)\in{\cal F}_{\rm reg}$ be a
fixed pair of polynomials, $V=V(f,g)$ the corresponding double
cover. Assume that the claim of Theorem 1 does not hold for $V$,
that is, for some effective divisor $D\in |-nK_V|$, $n\geq 1$, the
pair ($V,\frac{1}{n}D$) is not canonical. In the notations of Sec.
0.1, let $C=\varphi(E)\subset V$ be the centre of the singularity
$E$ on the variety $V$.

{\bf Proposition 0.1.} {\it The following inequality holds:
$\mathop{\rm codim}C\geq 3$}.

{\bf Proof.} Assume the converse: $C\subset V$ is a subvariety of
codimension 2. The Noether-Fano inequality implies immediately
that $\mathop{\rm mult}_CD> n$. In particular, the divisor $D$
does not coincide with the ramification divisor $R\subset V$ of
the double cover $\sigma$ (because $R$ is a smooth divisor on
$V$). Obviously, $R\cong W$, that is, $R$ is a smooth complete
intersection of codimension 2 in ${\mathbb P}$.

Restricting $D$ onto $R$, we obtain an effective divisor $Z$ on
$R$, which is cut out by a hypersurface of degree $n$ (in the
sense of the identification $R\cong W$, mentioned above). Let $Y$
be an irreducible component of the set $C\cap R$, $\mathop{\rm
codim}_RY=1$ (if $C\subset R$) or 2 (otherwise). We have the
inequality
$$
\mathop{\rm mult}\nolimits_YZ> n.
$$
Now the cone method (see [7, Proposition 3.6]) gives a
contradiction (the case when $\mathop{\rm codim}_RY=1$ is excluded
at once by the Lefschetz theorem for $R$), since $\mathop{\rm
dim}Y=M-3\geq 3$. Q.E.D. for the proposition.

Let $o\in C$ be a point of general position, $\varphi\,\colon
V^+\to V$ its blow up, $E=\varphi^{-1}(o)\subset V^+$ the
exceptional divisor, $D^+$ the strict transform of the divisor $D$
on $V^+$.

{\bf Proposition 0.2.} {\it For some hyperplane $B\subset E$ the
inequality
\begin{equation}\label{a1}
\mathop{\rm mult}\nolimits_oD+\mathop{\rm mult}\nolimits_BD^+>2n.
\end{equation}
holds.}

{\bf Proof} is given in [1, \S 3].

Therefore, the proof of Theorem 1 will be complete if we show that
for any point $o\in V$, where $V\in{\cal F}_{\rm reg}$ is a
regular variety, any (prime) divisor $D\in|-nK_V|$ and any
hyperplane $B\subset E$ the following inequality, which is
opposite to (\ref{a1}), holds:
\begin{equation}\label{a2}
\mathop{\rm mult}\nolimits_oD+\mathop{\rm mult}\nolimits_BD^+\leq
2n.
\end{equation}

A proof of this inequality makes the heart of this paper (\S\S
2-3).\vspace{0.3cm}

{\bf 0.5. Remarks.} The concept of birational (super)rigidity is
now universally accepted, but different authors use different
definitions [2,8-10]. It seems that today the most precise
definition is as follows: a projective rationally connected
variety $X$ is birationally rigid, if there exists a model
$\widetilde X$, which is birational to $X$, satisfying the
condition of Sec. 0.2, that is, for any movable system $\Sigma$ on
$\widetilde X$ there exists a birational self-map
$\chi\in\mathop{\rm Bir}\widetilde X=\mathop{\rm Bir}X$, providing
the equality of the thresholds $c_{\rm virt}(\Sigma)=c(\Sigma)$.
(As the example of [10] shows, on the very variety $X$ such a
self-map $\chi$ may not exist.) In other words, the variety $X$ is
birationally rigid, if on some its model $\widetilde X$ one can
untwist all maximal singularities of movable linear systems by
means of birational self-maps.

The present paper is based on the ideas and constructions
developed in [1]. In [4] we started to investigate Fano double
hypersurfaces, however the hardest case, when the point $o$ lies
on the ramification divisor of the morphism $\sigma$, was not
considered there. For this reason in [4] an intermediate result
was formulated on the canonicity of pairs ($V,\frac{1}{n}D)$,
where $V$ is a Fano double hypersurface.

Theorem 1 extends the action of the theorem on Fano direct
products [1] to the double hypersurfaces. It can be interpreted as
a statement on the global canonical threshold of the variety $V$:
for a generic double hypersurface $V$ the threshold is equal to 1.
The global log canonical thresholds are of importance in
differential geometry (existence of the K\" ahler-Einstein metric,
see [11-13]).

The claim of Theorem 1 for $m=1$ (the double spaces) was shown in
[1], for $m=2$ (the double quadrics) in [4].

Finally, combining the arguments of the present paper with the
constructions of [14] in the spirit of the paper [1], one can
prove in the same way the divisorial canonicity of iterated Fano
double covers
$$
V=V^{(l)}\to V^{(l-1)}\to\dots\to V^{(1)}\to Q\subset{\mathbb P},
$$
where $Q\subset{\mathbb P}^{M+k}$ is a Fano complete intersection
of index $\geq 3$, and each arrow is a double cover with a smooth
branch divisor, and moreover $V=V^{(l)}$ is a primitive Fano
variety.

%%%%%%%%%%%%%%%%%%%%%%%%%%%%%%%%%%%%%%%%%%%%%%%%%%%%%%%%%%%%%%%%
%%%%%%%%%%%%%%%%%%%%%%%%%%%%%%%%%%%%%%%%%%%%%%%%%%%%%%%%%%%%%%%%
%%%%%%%%%%%%%%%%%%%%%%% section 1

\section{Regular double covers}

In this section we formulate the regularity conditions, defining
the set ${\cal F}_{\rm reg}$. In \S 4 it is shown that it is
non-empty, that is, that the regularity conditions make sense. The
regularity conditions are local and differ essentially in the
cases when the given point $o\in V$ lies on the ramification
divisor of the morphism $\sigma$ and when it does
not.\vspace{0.3cm}

{\bf 1.1. The regularity conditions outside the branch divisor.}
We work in a fixed coordinate system $z_1,\dots,z_{M+1}$ on
${\mathbb P}$ with the origin at the point $p=\sigma(o)$. The
equation of the hypersurface $Q$ is $f=q_1+\dots+q_m$, the
equation of the hypersurface $W^*$, which cuts out on $Q$ the
branch divisor $W$, is $g=1+w_1+\dots+w_{2l}$. Following [3,14],
let us consider the formal series
$$
(1+t)^{1/2}= 1+\sum^{\infty}_{i=1}\gamma_i t^i= 1+\frac12
t-\frac18 t^2+\dots,
$$
where
$$
\gamma_i=(-1)^{i-1}\frac{(2i-3)!!}{2^ii!}=(-1)^{i-1}
\frac{(2i-3)!}{2^{2i-2}i!(i-2)!},
$$
and construct the corresponding formal series in the variables
$z_*$:
$$
\sqrt{g}= (1+w_1+\dots+w_{2l})^{1/2}=
1+\sum^{\infty}_{i=1}\gamma_i (w_1+\dots+w_{2l})^i=
$$
$$
=1+\sum^{\infty}_{i=1}\Phi_i(w_1,\dots,w_{2l}),
$$
where $\Phi_i(w_1(z_*),\dots,w_{2l}(z_*))$ are homogeneous
polynomials of degree $i$ in the variables $z_*$. Obviously,
$$
\Phi_i(w_*)=\frac12 w_i+ (\mbox{polynomial in}\,\,\,
w_{1},\dots,w_{i-1}).
$$
For instance, $\Phi_1(w_*)=\frac12 w_1$. Furthermore, for $j\geq
1$ set
$$
[\sqrt{g}]_j=1+\sum^j_{i=1}\Phi_i(w_*(z_*))
$$
and
$$
g^{(j)}=g-[\sqrt{g}]_j^2.
$$
It is easy to see that the first non-zero homogeneous component of
the polynomial $g^{(j)}$ is of degree $j+1$. Denote it by the
symbol $h_{j+1}$. Obviously,
\begin{equation}
\label{b1} h_{j+1}[g]=w_{j+1}+A_j(w_1,\dots,w_j),
\end{equation}
where the precise form of the polynomial $A_j$ is of no interest
for us. We say that the variety $V$ is {\it regular} at the point
$o$, if the following three regularity conditions (R1.1)-(R1.3)
hold.

(R1.1) The sequence of homogeneous polynomials
$$
q_1,\dots,q_m,h_{l+1},\dots h_{2l-1}
$$
is regular in the local ring ${\cal O}_{0,{\mathbb C}^{M+1}}$,
that is, the system of equations
$$
q_1=\dots=q_m=h_{l+1}=\dots=h_{2l-1}=0
$$
defines a one-dimensional set in ${\mathbb C}^{M+1}$, a finite set
of lines, passing through the origin.

(R1.2) The linear span of any irreducible component of the closed
algebraic set
$$
q_1=q_2=q_3=0
$$
in ${\mathbb C}^{M+1}$ is the hyperplane $q_1=0$.

(R1.3) Here we need to separate the two cases: $m=3$ and $m\geq
4$. For $m=3$ we require that the $\sigma$-preimage of a section
of any irreducible component of the closed algebraic set
$\overline{\{q_1=q_2=q_3=0\}}\subset Q$ by any anticanonical
divisor $\sigma^{-1}(P)\ni o$, containing the point $o$, should be
irreducible (that is, the corresponding double cover should not
break). Assume now that $m\geq 4$. In that case the regularity
condition requires that the closed algebraic set
$$
\sigma^{-1}(\overline{\{q_1=q_2=0\}}\cap Q)\subset V
$$
should be irreducible and any section of that set by an
anticanonical divisor $\sigma^{-1}(P)\ni o$, where $P\ni p$ is
some hyperplane,
\begin{itemize}
\item either is also irreducible and reduced,
\item or breaks into two irreducible components $B_1+B_2$, where
$B_i=\sigma^{-1}(Q\cap S_i)$ is the $\sigma$-preimage of a section
of $Q$ by a plane $S_i\subset{\mathbb P}$ of codimension three
and, moreover, $\mathop{\rm mult}_oB_i=\mathop{\rm mult}_pQ\cap
S_i=3$,
\item or non-reduced and is of the form $2B$, where
$B=\sigma^{-1}(Q\cap S)$ is the $\sigma$-preimage of a section of
$Q$ by a plane $S\subset{\mathbb P}$ of codimension three and,
moreover, $\mathop{\rm mult}_oB=\mathop{\rm mult}_p(Q\cap S)=3$.
\end{itemize}

{\bf Remark 1.1.} As we will see in the proof of Proposition 1.1
below, the second and third options in the condition (R1.3)
realize for a variety $V$ of general position only when the
quadric
$$
q_2\,|\,_{\{q_1=0\}\cap P}
$$
is respectively, a pair of planes $S_1\cup S_2$ or a double plane
$2S$. For this reason for $M\geq 7$ the second and third options
can be excluded.\vspace{0.3cm}

{\bf 1.2. The regularity conditions on the branch divisor.} Assume
that the point $o\in V$ lies on the ramification divisor of the
morphism $\sigma$, $p=\sigma(o)\in W$. Let $\varphi_V\colon V^+\to
V$ and $\varphi_Q\colon Q^+\to Q$ be the blow ups of the points
$o$ and $p$, respectively, $E_V=\varphi^{-1}_V(o)\subset V^+$ and
$E_Q=\varphi^{-1}_Q(p)\subset Q^+$ the exceptional divisors. Let
$W^+\subset Q^+$ be the strict transform of the hypersurface $W$,
$E_W=W^+\cap E_Q$ a hyperplane in $E_Q\cong{\mathbb P}^{M-1}$. The
symbols $E^*_V$, $E^*_Q$ and $E^*_W$ stand for the dual projective
spaces. The natural embedding $\sigma^*\colon
T^*_pW\hookrightarrow T^*_oV$ defines the embedding
$\sigma^*\colon E^*_W\hookrightarrow E^*_V$. The rational map
$\sigma^+\colon V^+\dashrightarrow Q^+$, induced by the morphism
$\sigma$, defines the linear projection $\sigma^+_E\colon
E_V\dashrightarrow E_W\cong{\mathbb P}^{M-2}$, which is dual to
that embedding.

In the coordinate form, let $z_1,\dots,z_{M+1}$ be affine
coordinates on ${\mathbb P}$ with the origin at the point $p$. The
hypersurface $Q$ is given by the equation
$$
f(z_1,\dots,z_{M+1})=q_1+\dots+q_m,
$$
where we can assume that $q_1\equiv z_{M+1}$, so that the
functions $z_1\,|\,_Q,\dots,z_M\,|\,_Q$ define on $Q$ a system of
local coordinates. The hypersurface $W^*$ is given by the equation
$$
g(z_1,\dots,z_{M+1})=w_1+\dots+w_{2l},
$$
where we may assume that $w_1=z_1$, so that the functions
$z_2\,|\,_W,\dots,z_M\,|\,_W$ define on the branch divisor $W$ a
system of local coordinates. The variety $V$ is given in ${\mathbb
A}^{M+2}_{(y,z_1,\dots,z_{M+1})}$ by the system of affine
equations
$$
y^2-g(z_*)=f(z_*)=0,
$$
which is locally of the form
$$
y^2-z_1+\dots=z_{M+1}+\dots=0,
$$
where the dots stand for the terms of order 2 and higher in $z_*$,
so that the functions $y\,|\,_V,z_2\,|\,_V,\dots,z_M\,|\,_V$ make
a system of local coordinates on $V$ with the origin at the point
$o$. In these coordinates the linear projection $\sigma^+_E$ takes
the form
$$
(y,z_2,\dots,z_M)\mapsto (z_2,\dots,z_M).
$$
Let $B\subset E_V$ be a hyperplane. There are two possible cases:

(A) either $B$ is given by an equation
$$
y+\lambda(z_2,\dots,z_M)=0,
$$
where $\lambda\in T^*_pW$ is some linear form (possibly a zero
one); in that case we say that $B$ {\it is not pulled back from}
$Q$,

(B) or $B$ is given by an equation
$$
\lambda(z_2,\dots,z_M)=0,
$$
where $\lambda\in T^*_pW\backslash\{0\}$ is some {\it non-zero}
linear form; in that case we say that $B$ {\it is pulled back
from} $Q$.

Obviously, $B$ is pulled back from $Q$ when and only when
$B\in\sigma^*(E^*_W)$ as a point of the dual projective space
$B\in E^*_V$. We say that the variety $V$ is {\it regular} at the
point $o$, where $p=\sigma(o)\in W$, if the following three
conditions (R2.1)-(R2.3) hold. Let us fix the system of
coordinates $z_1,\dots,z_{M+1}$ on ${\mathbb P}$ with the origin
at the point $p$, considered above. Restrictions of the
polynomials $w_i,q_j$ onto the plane $\{z_1=z_{M+1}=0\}$ we denote
by the symbols ${\bar w}_i,{\bar q}_j$.

(R2.1) For any linear form $\lambda(z_2,\dots,z_M)$ the sequence
\begin{equation}\label{b2}
\lambda^2(z_*)-{\bar w}_2,{\bar q}_2,\dots,{\bar q}_m
\end{equation}
is regular in ${\cal O}_{0,{\mathbb C}^{M-1}}$, that is the system
of equations
$$
\lambda^2(z_*)-{\bar w}_2={\bar q}_2=\dots={\bar q}_m=0
$$
defines a cone of codimension $m$ in ${\mathbb C}^{M-1}$ with the
vertex at the origin, and moreover, for the quadratic forms ${\bar
w}_2$ and ${\bar q}_2$ we have the estimates
$$
\mathop{\rm rk} {\bar w}_2\geq 4\,\, \mbox{and}\,\,\mathop{\rm rk}
{\bar q}_2 \geq 3.
$$

(R2.2) The linear span of any irreducible component of the closed
algebraic set
$$
{\bar q}_2={\bar q}_3=0
$$
in ${\mathbb C}^{M-1}_{(z_2,\dots,z_M)}$ is the whole space
${\mathbb C}^{M-1}$.

(R2.3) For $m=3$ the $\sigma$-preimage of any irreducible
component of the closed algebraic set
$$
\overline{\{\lambda(z_*)=z_1=z_{M+1}=q_2=q_3=0\}}\subset Q
$$
is irreducible for any linear form $\lambda(z_2,\dots,z_M)$. For
$m\geq 4$ for any linear form $\lambda(z_*)$ the closed set
$$
\sigma^{-1}\overline{(\{\lambda(z_*)=z_1=z_{M+1}=q_2=0\}}\cap Q
)\subset V
$$
\begin{itemize}
\item either is irreducible and reduced,

\item or breaks into two irreducible components
$B_1+B_2$, where $B_i=\sigma^{-1}(Q\cap S_i)$ is the
$\sigma$-preimage of a section of $Q$ by a plane
$S_i\subset\{\lambda(z_*)=z_1=z_{M+1}=0\}$ and, moreover,
$\mathop{\rm mult}_oB_i=3$,

\item or is non-reduced and is of the form $2B$, where
$B=\sigma^{-1}(Q\cap S)$ is the $\sigma$-preimage of a section of
$Q$ by a plane $S$ and, moreover, $\mathop{\rm mult}_oB=3$.
\end{itemize}
\vspace{0.3cm}

{\bf 1.3. Correctness of the regularity conditions.} It is obvious
from the explicit formulation of the regularity conditions, given
above, that the property of the variety $V$ to be regular at a
point $o$ depends on the pair of polynomials $(f,g)\in{\cal F}$
and the point $p=\sigma(o)\in{\mathbb P}$. Let ${\cal F}_{\rm
reg}(p)\subset{\cal F}$ be the set, consisting of pairs $(f,g)$
for which $f(p)=0$ and the double cover $V=V(f,g)$ is regular at
any (if there are two of them) point $o\in\sigma^{-1}(p)$. If
$g(p)\neq 0$, then the regularity is understood in the sense of
Sec. 1.1, whereas if $g(p)=0$, then in the sense of Sec. 1.2. Let
${\cal F}_{\rm reg}\subset{\cal F}$ be the set of pairs $(f,g)$,
for which the double cover $V=V(f,g)$ is regular at every point
$o\in V$.

{\bf Proposition 1.1.} {\it The set ${\cal F}_{\rm reg}$ contains
a non-empty Zariski open subset in ${\cal F}$.}

{\bf Proof} is given in \S 4. As usual (see the survey [2]), the
idea is to estimate the codimension of the closed set ${\cal
F}_{\rm non-reg}(p)$ of double covers, non-regular at the point
$p$. As soon as it is shown that $\mathop{\rm codim}_{\cal F}{\cal
F}_{\rm non-reg}(p)\geq M+2$, one can conclude that
$$
\overline{{\cal F}\backslash{\cal F}_{\rm
reg}}=\overline{\bigcup_{p\in{\mathbb P}}{\cal F}_{\rm
non-reg}(p)}
$$
is a proper closed subset of positive codimension in ${\cal F}$,
which immediately implies Proposition 1.1.

%%%%%%%%%%%%%%%%%%%%%%%%%%%%%%%%%%%%%%%%%%%%%%%%%%%%%%%%%%%%%%%%%%%%%%%%%%%%%%%%%
%%%%%%%%%%%%%%%%%%%%%%%%%%%%%%%%%%%%%%%%%%%%%%%%%%%%%%%%%%%%%%%%%%%%%%%%%%%%%%%%%
%%%%%%%%%%%%%%%%%%%%%%%%%%%%%%%%%%%%%%%%%%%Section 2

\section{Proof of inequality (\ref{a2}) for a point on the
ramification divisor}

In this section the inequality (\ref{a2}) is proved in the
assumption that the point $o\in V$ lies on the ramification
divisor, that is, $p=\sigma(o)\in W$. We use the notations of Sec.
1.2. Proof is obtained by two different methods, depending on
whether the hyperplane $B$ is pulled back from $Q$ or
not.\vspace{0.3cm}

{\bf 2.1. The hyperplane $B$ is not pulled back from $Q$.} In
order to prove the inequality (\ref{a2}), assume the converse: the
inequality (\ref{a1}) holds, where $D\in|-nK_V|$ is an irreducible
divisor, $D^+\subset V^+$ its strict transform. Let $\Lambda$ be a
pencil of hyperplanes in ${\mathbb P}$, generated by the tangent
hyperplanes $T_pQ$ and $T_pW^*$ (in the coordinate form it is the
pencil $\alpha z_1+\beta z_{M+1}=0$). Let $T\in\Lambda$ be a
general hyperplane.

Set $Q_T=Q\cap T$, $W_T=W\cap T$ and $V_T=\sigma^{-1}(Q_T)$.
Obviously, the hypersurface $Q_T$ is non-singular at the point
$p$, the divisor $W_T\subset Q_T$ has an isolated singularity at
the point $p$, so that $o\in V_T$ is an isolated singular point,
either, $\sigma_T=\sigma\,|\,_{V_T}\colon V_T\to Q_T$ is the
double cover, branched over $W_T$. With respect to the affine
coordinates $(y,z_2,\dots,z_{M+1})$ the variety $V_T$ is given by
the pair of equations
$$
y^2-w_2\,|\,_{\{z_1=-\frac{\beta}{\alpha}z_{M+1}\}}+\dots=z_{M+1}+
q_2\,|\,_{\{z_1= -\frac{\beta}{\alpha}z_{M+1}\}}+\dots,
$$
where the dots stand for the terms of order 3 and higher in
$z_2,\dots,z_{M+1}$. By the regularity condition (R2.1) the
quadric $w_2\,|\,_{\{z_1=z_{M+1}=0\}}$ has rank at least two, so
that $o\in V_T$ is an isolated quadratic singularity, where the
quadric $E_T=V^+_T\cap E_V$, where $V^+_T\subset V^+$ is the
strict transform of the divisor $V_T$, is given by the equation
$$
y^2-w_2\,|\,_{\{z_1=z_{M+1}=0\}}=0.
$$
The latter quadratic form is of rank at least three, so that $E_T$
does not contain hyperplanes in $E_V$. Therefore, $\mathop{\rm
mult}_oV_T+\mathop{\rm mult}_BV^+_T=2$. Since $V_T\in|-K_V|$ is an
anticanonical divisor, we get $D\neq V_T$. Thus the effective
cycle $D_T=(D\circ V_T)$ of codimension two is well defined, it is
an effective divisor on $V_T$. Obviously, $\mathop{\rm
mult}_oD_T\geq 2\mathop{\rm mult}_oD$, and the strict transform
$D^+_T\subset V^+$ contains the hyperplane section $B_T=B\cap
V^+_T$ of the quadric. The quadric $B_T$ is given on $E_T$ by the
equation $y+\lambda(z_2,\dots,z_M)=0$. By the regularity condition
(R2.1), the quadric $B_T$ is irreducible and reduced.

{\bf Lemma 2.1.} {\it The following inequality holds:}
\begin{equation}\label{c1}
\mathop{\rm mult}\nolimits_oD_T+2\mathop{\rm
mult}\nolimits_{B_T}D^+_T>4n.
\end{equation}

{\bf Proof.} One has to prove this inequality, because the tangent
cone of the divisor $D$, that is, the algebraic cycle ($D^+\circ
E_V$), can contain the quadric $E_T$. By the standard formulas of
the intersection theory [15],
\begin{equation}\label{c2}
\mathop{\rm mult}\nolimits_oD_T=2\mathop{\rm mult}\nolimits_oD+2a,
\end{equation}
where the integer $a\in{\mathbb Z}_+$ is defined by the relation
$$
(D^+\circ V^+_T)=aE_T+Z,
$$
the effective cycle $Z$ does not contain $E_T$. On the other hand,
every irreducible component of the cycle $Z$ is not contained in
the exceptional divisor $E_V$, that is, $Z=(D\circ V_T)^+$ and
thus
\begin{equation}\label{c3}
\mathop{\rm mult}\nolimits_{B_T}D^+_T+a\geq\mathop{\rm
mult}\nolimits_BD^+.
\end{equation}
Combining the estimates (\ref{c2}) and (\ref{c3}), we get
(\ref{c1}). Q.E.D. for the lemma.

Now let us consider the standard hypertangent linear systems
$$
\Lambda^Q_1,\dots,\Lambda^Q_{m-1}
$$
on the hypersurface $Q$ at the point $p$. In terms of the affine
coordinates $z_*$ these systems take the form
$$
\Lambda^Q_i=\left|\sum^i_{j=1}s_{i-j}(q_1+q_2+\dots+q_j)=0\right|,
$$
where $s_a(z_*)$ are arbitrary homogeneous polynomials of degree
$a\in{\mathbb Z}_+$ in $z_*$, $f_j=q_1+\dots+q_j$ is the left
segment of length $j$ of the equation $f$. Set
$\Lambda_i=\sigma^*\Lambda^Q_i$ to be the pull back of these
systems on $V$, $\Lambda^+_i$ their strict transform on $V^+$.
Obviously,
$$
\Lambda_i\subset|-iK_V|,\quad \mathop{\rm
mult}\nolimits_o\Lambda_i\geq i+1,
$$
and by the regularity condition (R2.1) we can say more precisely
that $\mathop{\rm mult}_o\Lambda_i=i+1$, that is, $\Lambda^+_i$ is
a subsystem of the complete linear system
$|-i\varphi^*_VK_V-(i+1)E_V|$. Set also
$\Lambda^E_i=\Lambda^+_i\,|\,_{E_V}$. These are linear systems of
hypersurfaces of degree $(i+1)$ on $E_V\cong{\mathbb P}^{M-1}$.

The regularity condition (R2.1) implies that the base sets of the
linear systems $\Lambda_i$ and $\Lambda^E_i$ are of codimension
$i$ with respect to $V$ and $E_V$, respectively, and moreover, the
codimension does not change when we restrict on $V_T$ and $E_T$,
respectively. These base sets are given by the system of equations
$$
q_2|\,_{\{z_1=z_{M+1}=0\}}=\dots=q_{i+1}|\,_{\{z_1=z_{M+1}=0\}}=0.
$$
Let
$(R_2,\dots,R_{m-1})\in\Lambda_2\times\dots\times\Lambda_{m-1}$ be
a generic set of divisors of the hypertangent systems. Their
strict transforms on $V^+$ and restrictions onto $E_V$ are denoted
by the symbols $R^+_i$ and $R^E_i=(R^+_i\circ E_V)$, respectively.
By the regularity condition the set-theoretic intersections
$$
\mathop{\rm Supp}D_T\cap R_2\cap\dots\cap R_i
$$
and
$$
\mathop{\rm Supp}D^+_T\cap E_V\cap R^E_2\cap\dots\cap R^E_i
$$
are on codimension $i$ in $V_T$ and $E_T$, respectively,
$i=2,\dots,m-1$. Therefore, the effective cycles
$$
Y_{i+1}=(D_T\circ R_2\circ\dots\circ R_i)
$$
and
$$
Y^E_{i+1}=(D^+_T\circ E_V\circ R^E_2\circ\dots\circ R^E_i)
$$
are well defined, and moreover, $Y^E_i=(Y^+_i\circ E_V)$, that is,
$Y^E_i$ is the projectivized tangent cone to $Y_i$ at the point
$o$. In particular,
$$
\mathop{\rm mult}\nolimits_oY_i=\mathop{\rm
deg}Y^E_i=\frac{i!}{2}\mathop{\rm mult}\nolimits_oD_T,
$$
where the degree of the cycle $Y^E_i$ is understood in the sense
of the projective space $E_V\cong{\mathbb P}^{M-1}$.\vspace{0.3cm}

{\bf 2.2. End of the proof.} Let us use the method of [5] (in a
modified form).

{\bf Definition 2.1.} An irreducible component $Z$ of the cycle
$Y_i$ is a $B$-{\it component}, if at least one irreducible
component of its projectivized tangent cone at the point $o$, that
is, of the algebraic cycle $(Z^+\circ E_V)$, where $Z^+$ is the
strict transform of $Z$ on $V^+$, is contained in the hyperplane
$B$.

For the effective cycle $Y_m$ we have the decomposition
$$
Y_m=Y_{\sharp}+Y_B,
$$
where in $Y_B$ we collect all $B$-components of the cycle $Y_m$,
and only them. The following fact is crucial.

{\bf Lemma 2.2.} {\it None of the irreducible components of the
cycle $Y_B$ is contained in the inverse image of the tangent
hyperplane $\sigma^{-1}(T_pQ\cap Q)=\{z_{M+1}=0\}$.}

{\bf Proof.} Assume the converse: such a $B$-component $Z$ does
exist, $\sigma(Z)\subset T_pQ$. Then its tangent cone
$Z^E=(Z^+\circ E_V)$ is entirely contained in the tangent cone of
the divisor $\sigma^{-1}(T_pQ\cap Q)$. The latter is given in
$E_V$ by the equation
$$
q_2|\,_{\{z_1=z_{M+1}=0\}}=0.
$$
By the definition of a $B$-component, there is an irreducible
component $S$ of the cycle $Z^E$, which is contained in $B$.
Therefore, the polynomials
$$
y+\lambda(z_2,\dots,z_M),y^2-{\bar w}_2, {\bar q}_2,\dots,{\bar
q}_m
$$
vanish on $S$. Recall that $S\subset E_V$ and $(y:z_2:\dots:z_M)$
are homogeneous coordinates on the projective space $E_V$.
Therefore, on $S$ vanish the polynomials
\begin{equation}\label{c4}
\lambda^2(z_2,\dots,z_M)-{\bar w}_2, {\bar q}_2,\dots,{\bar q}_m.
\end{equation}
These polynomials do not depend on $y$, so that the set of their
common zeros is a set of cones in $E_V\cong{\mathbb P}^{M-1}$ with
the vertex at the point $(1:0\dots:0)$. On the other hand, the
hyperplane $B$ does not contain that point, so that $S$ does not
contain it, too. Therefore, $m$ homogeneous polynomials (\ref{c4})
in the variables $z_2,\dots,z_M$ vanish on the irreducible
subvariety $\bar{S}\subset{\mathbb P}^{M-2}$, the projection of
$S$ from the point $(1:0:\dots:0)$. By what was said above,
$\mathop{\rm dim}\bar{S}=\mathop{\rm dim}S$, so that $\mathop{\rm
codim}\bar{S}=m-1$, which contradicts the regularity condition
(R2.1). Q.E.D. for Lemma 2.2.

Let us come back to the proof of inequality (\ref{a2}). By
construction, the algebraic cycle $Y_m$ is of degree
$$
\mathop{\rm deg}Y_m=\mathop{\rm deg}Y_{\sharp}+\mathop{\rm
deg}Y_B=n\cdot 2m\cdot (m-1)!=2nm!
$$
and its multiplicity at the point $o$ equals to
\begin{equation}\label{c5}
\mathop{\rm mult}\nolimits_oY_m=\mathop{\rm mult}\nolimits_oY_{
\sharp}+\mathop{\rm mult}\nolimits_oY_B=\frac{m!}{2}\mathop{\rm
mult}\nolimits_oD_T.
\end{equation}
For the cycle $Y_{\sharp}$ we have the universal estimate
$\mathop{\rm mult}_o Y_{\sharp}\leq\mathop{\rm deg}Y_{\sharp}$,
which cannot be improved. The situation with the cycle $Y_B$ is
much better: by Lemma 2.2, the effective cycle
$$
Y^*=(Y_B\circ\sigma^{-1}(T_pQ\cap Q))
$$
is well defined, and for that cycle we get
$$
\mathop{\rm deg}Y^*=\mathop{\rm deg}Y_B,\,\quad\mathop{\rm
mult}\nolimits_oY^*\geq 2\mathop{\rm mult}\nolimits_oY_B,
$$
so that by the universal inequality $\mathop{\rm
mult}_oY^*\leq\mathop{\rm deg}Y^*$ we get the estimate
\begin{equation}\label{c6}
\mathop{\rm mult}\nolimits_oY_B\leq\frac12\mathop{\rm deg}Y_B.
\end{equation}
Finally, by the construction of the cycle $Y_m$, we have the
inequality
$$
\mathop{\rm mult}\nolimits_oY_B=\mathop{\rm
deg}Y^E_B\geq\frac{m!}{2}\mathop{\rm deg}B_T\mathop{\rm
mult}\nolimits_{B_T}D^+_T=m!\mathop{\rm mult}\nolimits_{B_T}D^+_T.
$$
Combining this inequality with the estimates (\ref{c5}) and
(\ref{c6}), we get finally
$$
(2\mathop{\rm mult}\nolimits_{B_T}D^+_T)m!+(\mathop{\rm
mult}\nolimits_oD_T)m!\leq
$$
$$
\leq\mathop{\rm deg}Y_B+2\mathop{\rm
mult}\nolimits_oY_{\sharp}+2\mathop{\rm mult}\nolimits_oY_B\leq
$$
$$
\leq\mathop{\rm deg}Y_B+2\mathop{\rm deg}Y_{\sharp}+\mathop{\rm
deg}Y_B=2\mathop{\rm deg}Y_m=4nm!,
$$
whence, reducing by $m!$, we get the inequality, which is opposite
to (\ref{c1}). This proves the inequality (\ref{a2}) in the case
when the hyperplane $B\subset E_V$ is not pulled back from
$Q$.\vspace{0.3cm}

{\bf 2.3. The hyperplane $B$ is pulled back from $Q$.} Let us
prove the inequality (\ref{a2}), assuming that the hyperplane
$B\subset E_V$ is pulled back from $Q$, that is, in the
coordinates $(y,z_2,\dots,z_M)$ it has an equation
$\lambda(z_2,\dots,z_M)=0$, where $\lambda\not\equiv 0$ is some
linear form. Again let us assume the converse: $\mathop{\rm
mult}_oD+\mathop{\rm mult}_BD^+>2n$, where $D\in|-nK_V|$ is an
irreducible divisor, $D^+\subset V^+$ its strict transform. Here
we argue following the model of the case when the point $o\in V$
lies outside the ramification divisor.

Let $\Lambda_B$ be the two-dimensional system of hyperplanes in
${\mathbb P}$, cutting out on the tangent space $T_pW$ the
hyperplane $B$. In the coordinates $(z_1,\dots,z_{M+1})$ the
equations of these hyperplanes are of the form
$$
\alpha z_1+\beta z_{M+1}+\lambda(z_2,\dots,z_M)=0,
$$
$\alpha,\beta\in{\mathbb C}$ are constants. Let
$\Lambda=\sigma^*(\Lambda_B\,|\,_Q)$ be the corresponding linear
system on $V$, $R\in\Lambda$ a general divisor, $Q_R=\sigma(R)$ a
smooth hypersurface of degree $m$ in ${\mathbb P}^M$. The branch
divisor $W_R=Q_R\cap W$ and the variety $R$ itself are smooth at
the points $p$ and $o$, respectively.

{\bf Lemma 2.3.} {\it For the divisor $D_R=D\cap R$ the estimate
$\mathop{\rm mult}_oD_R>2n$ holds.}

{\bf Proof.} By the formulas of the elementary intersection theory
[15] we get: $\mathop{\rm mult}_oD_R=\mathop{\rm
mult}_oD+\mathop{\rm deg}Z$, where the effective divisor $Z$ on
$E_V$ is defined by the relation $(D^+\circ R^+)=D^+_R+Z$. By
construction of the divisor $R$, we get $(R^+\circ E_V)=B$, so
that $Z$ contains $B$ with multiplicity at least $\mathop{\rm
mult}_BD^+$. This proves the lemma.

Consider now the pencil $\Lambda_R$ of hyperplane sections of
$Q_R$, tangent to the branch divisor $W_R$ at the point $p$. In
the coordinate form it is the pencil of hyperplanes
$\alpha^*z_1+\beta^*z_{M+1}=0$, restricted onto $Q_R$. Let
$T\in\Lambda_R$ be a general divisor of that pencil. Set
$W_T=W\cap T$ and $V_T=\sigma^{-1}(T)$. The divisor $T$ is
non-singular at the point $p$, whereas the divisor $W_T$ has at
that point an isolated quadratic singularity.

More precisely, the functions $z_2,\dots,z_M$ form a system of
local coordinates on $Q_R$, the tangent hyperplane to the branch
divisor $W_R$ is given by the equation $\lambda(z_2,\dots,z_M)=0$,
so that the linear component of the local equation of the divisor
$T$ in these coordinates is $\lambda(z_2,\dots,z_M)$. The tangent
cone to the divisor $V_T$ in the coordinates $(y:z_2:\dots:z_M)$
on $E_V$ is given by the pair of equations
$$
y^2-w_2|\,_{\{z_1=z_{M+1}=0\}}=\lambda(z_*)=0.
$$
Therefore, $\mathop{\rm mult}_oV_T=2$. The divisor $V_T$ is
irreducible and for this reason $D_R\neq V_T$ by Lemma 2.3.
Therefore, the effective cycle
$$
D_T=(D_R\circ V_T),
$$
satisfying the estimate $\mathop{\rm mult}_oD_T>4n$, is well
defined. The effective divisor $D_T$ on the singular double cover
$V_T$ can be assumed to be irreducible.

{\bf Lemma 2.4.} {\it The divisor
$$
S=\sigma^{-1}(T_pQ\cap T)=\sigma^{-1}(T_pT)
$$
on the variety $V_T$ is irreducible and has multiplicity precisely
4 at the point $o$.}

{\bf Proof.} By the regularity condition (R2.3) the system of
three equations
$$
y^2-w_2|\,_{\{z_1=z_{M+1}=0\}}=\lambda(z_*)=q_2|\,_{\{z_1=z_{M+1}=0\}}=0
$$
defines on $E_V\cong{\mathbb P}^{M-1}_{(y:z_2:\dots:z_M)}$ an
effective cycle of codimension 3 and degree 4. Q.E.D. for the
lemma.

From the lemma that we have just proven, it follows that
$D_T\not\subset\sigma^{-1}(T_pQ)$, or, equivalently, $D_T\neq S$.
Consider the second hypertangent system
$$
\Lambda_2=|s_0f_2+s_1f_1|,
$$
where recall that $f_1=q_1$, $f_2=q_1+q_2$, $s_i(z_*)$ are
homogeneous polynomials in the variables $z_1,\dots,z_{M+1}$. By
the regularity condition (R2.3), the base set of its restriction
$\Lambda^T_2=\Lambda_2|\,_T$ has codimension two. Thus for a
general divisor $L\in\Lambda^T_2$ we get $\sigma(D_T)\not\subset
L$, so that the effective cycle of codimension two on $V_T$
$$
D_L=(D_T\circ\sigma^{-1}(L))
$$
is well defined and satisfies the estimate $\mathop{\rm
mult}_oD_L>12n$. The degree of the cycle $D_L$ (in the sense of
the anticanonical class $(-K_V)$) is $2n\mathop{\rm deg}V=4nm$.
Thus
\begin{equation}\label{c7}
\frac{\mathop{\rm mult}_o}{\mathop{\rm
deg}}D_L>\frac{6}{\mathop{\rm deg}V}.
\end{equation}
However, it follows from the regularity condition (R2.3), that the
base set of the linear system $\sigma^*\Lambda^T_2$ either is
irreducible and satisfies the equality
$$
\frac{\mathop{\rm mult}_o}{\mathop{\rm deg}}\mathop{\rm
Bs}\sigma^*\Lambda^T_2=\frac{6}{\mathop{\rm deg}V},
$$
or breaks into two components with the same ratio $\mathop{\rm
mult}_o/\mathop{\rm deg}$. Replacing the cycle $D_L$ by its
suitable irreducible component, we may assume that $D_L$ is
irreducible and by the inequality (\ref{c7})
$D_L\not\subset\mathop{\rm Bs}\sigma^*\Lambda^T_2$. By the
construction of the linear system $\Lambda_2$, some polynomial of
the form $u_0f_2+u_1f_1$ vanishes on $\sigma(D_L)$, where $u_0\neq
0$ is a constant, $u_1(z_*)$ is a linear form. Without loss of
generality assume that $u_0=1$. It follows that
\begin{equation}\label{c8}
f_2|\,_{\sigma(D)}\equiv-u_1q_1|\,_{\sigma(D)}.
\end{equation}

{\bf Lemma 2.5.} {\it The image $\sigma(D_L)$ is not contained in
$T_pQ$.}

{\bf Proof.} The claim of the lemma means that $q_1=f_1=z_{M+1}$
does not vanish on $\sigma(D_L)$. Assume the converse. Then
(\ref{c8}) implies that $\sigma(D_L)\subset\mathop{\rm
Bs}\Lambda^T_2$, but we already know that this is not the case.
Q.E.D. for the lemma.

Consider the effective cycle
$$
D^{\sharp}=(D_L\circ\sigma^{-1}(T_pQ)).
$$
It is of codimension 4 on $V$ and satisfies the inequality
$$
\frac{\mathop{\rm mult}_o}{\mathop{\rm
deg}}D^{\sharp}>\frac{12}{\mathop{\rm deg}V}.
$$
We may assume that the cycle $D^{\sharp}$ is an irreducible
subvariety in $V$. Its image $\Delta=\sigma(D^{\sharp})$ is a
subvariety of codimension 4 on $Q$, satisfying the estimate
$$
\frac{\mathop{\rm mult}_p}{\mathop{\rm
deg}}\Delta>\frac{6}{\mathop{\rm deg}Q}.
$$
The hypersurface $Q$ satisfies the regularity condition (R2.1), so
that, repeating the arguments of [16, Sec. 4] word for word (see
also [2]), we get a contradiction. The proof of the inequality
(\ref{a2}) in the case when $p\in W$ and the hyperplane $B\subset
E_V$ is pulled back from $Q$, is complete.

Note that the last contradiction, completing the proof, can be
obtained directly on $V$, without pushing $D^{\sharp}$ down on
$Q$, but, on the contrary, pulling back the hypertangent systems
$$
\Lambda^Q_i=\left|\sum^i_{j=1}s_{i-j}f_j=0\right|
$$
on $V$, as it was done above for the second hypertangent system.
%%%%%%%%%%%%%%%%%%%%%%%%%%%%%%%%%%%%%%%%%%%%%%%%%%%%%%%%%%%%%%%%%%%%%%%%%%%%%%%%%%%%
%%%%%%%%%%%%%%%%%%%%%%%%%%%%%%%%%%%%%%%%%%%%%%%%%%%%%%%%%%%%%%%%%%%%%%%%%%%%%%%%%%%%
%%%%%%%%%%%%%%%%%%%%%%%%%%%%%%%%%%%%%%%%%%Section 3
\section{Proof of the inequality (\ref{a2}) for a point outside
the ramification divisor}

We use the notations and conventions of Sec. 1.1. Let us prove the
inequality (\ref{a2}), assuming that the point $o\in V$ lies
outside the ramification divisor of the morphism $\sigma$, that
is, $p=\sigma(o)\in Q$ lies outside the divisor $W$. In that case
we argue precisely following the scheme of the Fano hypersurfaces
[1, \S 2, Sec. 2.1]. Since our constructions are almost word for
word the same as those in [1], we give only the principal steps of
the proof.

Let $\varphi\colon V^+\to V$ and $\varphi_Q\colon Q^+\to Q$ be the
blow ups of the points $o$ and $p$ on $V$ and $Q$, respectively,
$E\subset V^+$ and $E_Q\subset Q^+$ the exceptional divisors. The
morphism $\sigma$ extends to a regular map
$V^+\backslash\{o_1\}\to Q^+$, identifying $E$ and $E_Q$, where
$\sigma^{-1}(p)=\{o,o_1\}$. This identification $E\cong
E_Q\cong{\mathbb P}^{M-1}$ will be meant in the sequel without
special reservations. Assume that the inequality $\mathop{\rm
mult}_oD+\mathop{\rm mult}_BD^+>2n$ holds, where $B\subset E$ is a
hyperplane, $D\in|-nK_V|$, $D^+\subset V^+$ is the strict
transform. The divisor $D$ is assumed to be irreducible and
reduced. Let $B_Q\subset E_Q$ be the corresponding hyperplane in
$E_Q$. The exceptional divisor $E_Q$ identifies naturally with the
projectivization ${\mathbb P}(T_pQ)$. Let $\Lambda_B$ be the
pencil of hyperplanes in ${\mathbb P}$, cutting out $B$ on $E_Q$
and $\Lambda=\sigma^*(\Lambda_B|\,_Q)$ the pull back on $V$ of its
restriction onto $Q$. Consider a general divisor $R\in\Lambda$,
let $R^+\subset V^+$ be its strict transform. The divisor $R$ is
smooth at the point $o$, and moreover,
$$
R^+\cap E=B.
$$
Set $D_R=D|\,_R=(D\circ R)$. This is an effective divisor on the
variety $R$.

{\bf Lemma 3.1.} {\it The following inequality holds:}
\begin{equation}\label{d1}
\mathop{\rm mult}\nolimits_oD_R>2n.
\end{equation}

{\bf Proof} is word for word the same as the proof of Lemma 3 in
[1, \S 2, Sec. 2.1].

{\bf Lemma 3.2.} {\it The divisor
$T_R=\sigma^{-1}(T_p\sigma(R))\cap R$ on the variety $R$ is
irreducible and has multiplicity precisely 2 at the point $o$.}

{\bf Proof} is word for word the same as the proof of Lemma 4 in
[1, \S 2, Sec. 2.1], based on the regularity condition (R1.2).

By Lemmas 3.1 and 3.2, we may assume the divisor $D_R$ to be
irreducible and reduced, and different from $T_R$. Consider the
second hypertangent system on $Q$:
$$
\Lambda^Q_2=|s_0f_2+s_1f_1=0|\,_Q,
$$
$s_i(z_*)$ are homogeneous polynomials of degree $i$. The base set
of its restriction $\Lambda^R_2=\sigma^*\Lambda^Q_2\,|\,_R$ onto
$R$, that is,
$$
S_R=\{\sigma^*q_1|\,_R=\sigma^*q_2|\,_R=0\},
$$
has by the regularity condition (R1.3) codimension 2 in $R$ and
either is irreducible and of multiplicity 6 at the point $o$, or
breaks into two $\sigma$-preimages of plane sections of the
hypersurface $\sigma(R)$, each of multiplicity 3 at the point $o$.
In any case, for a general divisor $L\in\Lambda^R_2$ we get
$D_R\not\subset L$, so that the effective cycle of codimension two
$D_L=(D_R\circ L)$ is well defined and satisfies the estimate
\begin{equation}\label{d2}
\frac{\mathop{\rm mult}_o}{\mathop{\rm
deg}}D_L>\frac{3}{\mathop{\rm deg}V}=\frac{3}{2m}.
\end{equation}
Replacing the cycle $D_L$ by its suitable irreducible component,
one may assume it to be an irreducible subvariety of codimension 2
in $R$, and comparing the estimate (\ref{d2}) with the description
of the set $S_R$ given above, we see that $D_L\not\subset S_R$.
Now the arguments of Sec. 2.1 in [1] show that this implies that
$D_L\not\subset T_R$ (similar arguments are used in the present
paper in Sec. 2.3, see the proof of Lemma 2.5). It follows that
the effective cycle
$$
Y=(D_L\circ T_R)
$$
of codimension 4 on $V$ is well defined and satisfies the estimate
\begin{equation}\label{d3}
\frac{\mathop{\rm mult}_o}{\mathop{\rm deg}}Y>\frac{6}{\mathop{\rm
deg}V}=\frac{3}{m}.
\end{equation}

Now a contradiction is achieved by the method of the paper [3]. By
linearity of the multiplicity and degree, the cycle $Y$ can be
assumed to be an irreducible variety. In the notations of Sec. 1.1
let
\begin{equation}\label{d4}
\Lambda_k=\left|\sum^K_{i=1}s_{k-i}f_i+
\sum^k_{i=l}s^*_{k-i}(y-[\sqrt{g}])_i\right|,
\end{equation}
$k=1,2,\dots $, be the $k$-th hypertangent system, where
$s_j,s^*_j$ are homogeneous polynomials in $z_*$ of degree $j$,
for simplicity of notations we omit the symbol $\sigma^*$
(strictly speaking, we should have written
$\Sigma\sigma^*s_{k-i}\sigma^*f_i$ etc.) and, finally, the
right-hand sum in (\ref{d4}) is assumed to be equal to zero, if
$k<l$. The system $\Lambda_k$ is a subsystem of the complete
system $|-kK_V|$ and it is easy to see that $\mathop{\rm
mult}_o\Lambda_k\geq k+1$.

Set
$$
{\cal M}=[1,m-1]\cap{\mathbb Z}_+=\{1,\dots,m-1\}
$$
and ${\cal L}=[l,2l-2]\cap{\mathbb Z}_+ =\{l,\dots,2l-2\}$. By the
regularity condition, for the codimension of the base set of the
hypertangent system $\Lambda_k$ we get
$$
\mathop{\rm codim}\mathop{\rm Bs}\Lambda_k= \sharp [1,k]\cap{\cal
L}+\sharp [1,k]\cap{\cal M}.
$$
Let
$$
(D_1,\dots,D_{m-1},D^*_l,\dots,D^*_{2l-2})\in \prod_{k\in{\cal
M}}\Lambda_k\times\prod_{k\in{\cal L}}\Lambda_k
$$
be a general set of hypertangent divisors. Re-order this set as
$(L_1,\dots,L_{M-1})$ in such a way that for all divisors
$L_i\in\Lambda_{k(i)}$ the inequality $k(i+1)\geq k(i)$ holds. Now
the regularity condition implies that
$$
\mathop{\rm codim}\nolimits_oY\cap L_5\cap\dots\cap L_k=k,
$$
where $\mathop{\rm codim}_o$ means the codimension in a
neighborhood of the point $o$ with respect to $V$. Therefore, one
can realize the standard procedure of constructing irreducible
subvarieties $Y=Y_4,Y_5,\dots,Y_{M-1}$ of codimension $\mathop{\rm
codim}Y_i=i$, where $Y_{i+1}$ is an irreducible component of the
effective cycle $(Y_i\circ L_{i+1})$ with the maximal value of the
ratio $\mathop{\rm mult}_o/\mathop{\rm deg}$. The subvariety
$Y^*=Y_{M-1}$ is an irreducible curve on $V$, satisfying the
inequality
\begin{equation}\label{d5}
\frac{\mathop{\rm mult}_o}{\mathop{\rm
deg}}Y^*\geq\frac{\mathop{\rm mult}_o}{\mathop{\rm deg}}Y\cdot
\left(\prod^{M-1}_{i=5}\frac{k(i)+1}{k(i)}\right).
\end{equation}
It is not hard to check that the value of the product in the
brackets is
$$
\frac{m}{5}\cdot\frac{2l-1}{l}=\frac{m(2l-1)}{5l}>\frac{m}{3}
$$
for $l\geq 4$, $3m/8$ for $l=3$ and $m/3$ for $l=2$. In any case
this value is not less than $m/3$. Since the left-hand side of the
inequality (\ref{d5}) is not higher than one, we obtain for the
ratio $(\mathop{\rm mult_o}/\mathop{\rm deg})Y$ the estimate,
which is opposite to the inequality (\ref{d3}). This contradiction
completes the proof of inequality (\ref{a2}).
%%%%%%%%%%%%%%%%%%%%%%%%%%%%%%%%%%%%%%%%%%%%%%%%%%%%%%%%%%%%%%%%%%%%%%%%%%%%%%%%%%%
%%%%%%%%%%%%%%%%%%%%%%%%%%%%%%%%%%%%%%%%%%%%%%%%%%%%%%%%%%%%%%%%%%%%%%%%%%%%%%%%%%%
%%%%%%%%%%%%%%%%%%%%%%%%%%%%%%%%%%%Section 4

\section{Proof of the regularity conditions}

In this section, we prove Proposition 1.1.\vspace{0.3cm}

{\bf 4.1. The regularity conditions outside the branch divisor.}
The condition (R1.1) for a generic cover $V$ was shown in [3]. To
prove the condition (R1.2) for a generic hypersurface $Q$, one
needs to argue word for word in the same way as for this condition
for a generic Fano hypersurface, see [1,\S 2], because in this
condition only the three components $q_1,q_2,q_3$ take part.
However, the condition (R1.3) for a double cover is essentially
different from the corresponding condition for a hypersurface and
for this reason needs a special consideration.

Assume at first that $m=3$. Consider the following general
situation: $X\subset{\mathbb P}^N$ is an irreducible subvariety,
$x\in X$ a smooth point, $\mathop{\rm dim}X=k\geq 2$.

{\bf Lemma 4.1.} {\it The closed subset $\Xi(x)\subset{\cal
P}_{2l}(x)$ in the space of homogeneous polynomials of degree $2l$
in the homogeneous variables on ${\mathbb P}^N$, vanishing at the
point $x$, defined by the condition
$$
g\in \Xi(x)\Leftrightarrow g=h^2\,\,\mbox{in}\,\,{\cal O}_{x,X}
$$
is of codimension at least ${2l+k-1}\choose {k-1}$}.

{\bf Proof.} Let $(u_1,\dots,u_N)$ be some system of affine
coordinates on ${\mathbb P}^N$ with the origin at the point $x$,
whereas $(u_1,\dots,u_k)$ make a system of local parameters on $X$
at that point. Consider the standard projection
\begin{equation}\label{e1}
\pi\colon{\cal O}_{x,X}\to{\cal O}_{x,X}/ {\cal
M}^{2l+1}_x\cong{\mathbb C}[u_1,\dots,u_k]/
(u_1,\dots,u_k)^{2l+1},
\end{equation}
where ${\cal M}_x=(u_1,\dots,u_k)\subset{\cal O}_{x,X}$ is the
maximal ideal of the local ring. We denote the latter algebra in
(\ref{e1}) by the symbol ${\cal A}_{2l}$, and its maximal ideal by
the symbol ${\cal M}$. Let $\Xi\subset{\cal M}$ be the set of full
squares,
$$
\Xi=\{g\in{\cal M}|\,g=h^2\,\,\mbox{for some}\,\,h\in{\cal M}\}.
$$
Obviously, $\Xi(x)\subset\pi^{-1}(\Xi)$. Furthermore, since the
restriction of $\pi$ onto ${\cal P}_{2l}(x)$ is a linear
surjective map, we get the estimate
$$
\mathop{\rm codim}\nolimits_{{\cal P}_{2l}(x)}\Xi(x)\geq
\mathop{\rm codim}\nolimits_{{\cal M}}\Xi.
$$
It is not hard to estimate the latter codimension from below. Let
$$
h=h_1+\dots+h_{2l}\in{\cal M}
$$
be an arbitrary element, decomposed into homogeneous components,
that is, homogeneous polynomials in $u_1,\dots,u_k$. For its
square we have the presentation
$$
h^2=\sum^{2l}_{i=2}v_i(h_1,\dots,h_{i-1}),
$$
where the homogeneous component $v_i$ of degree $i$ depends on
$h_1,\dots,h_{i-1}$ only. In particular, $h^2$ does not depend on
the last component $h_{2l}$ and for this reason
$$
\mathop{\rm dim}\Xi\leq\mathop{\rm dim}{\cal M}/{\cal M}^{2l}.
$$
Therefore, $\mathop{\rm codim}_{{\cal M}}\Xi$ is not less than the
dimension of the space ${\cal M}^{2l}/{\cal M}^{2l+1}$, which is
equal to ${2l+k-1}\choose {k-1}$. Q.E.D. for the lemma.

Let us come back to the regularity condition (R1.3) for the double
cubic ($m=3$). Let $Y\subset{\mathbb P}$ be an irreducible
component of the closed set
$$
\overline{\{q_1=q_2=q_3=0\}}\cap P,
$$
where $P\ni p$ is a hyperplane. Obviously, $Y$ is a cone with the
vertex at the point $p$. The hypersurface $W^*_{2l}$, that cuts
out on $Q$ the branch divisor, does not contain the point $p$.
Therefore,
$$
W^*\cap Y\not\subset\mathop{\rm Sing}Y
$$
and we can apply Lemma 4.1 to $X=Y$. We obtain that a violation of
the regularity condition (R1.3) at the point $p$ imposes on the
polynomial $g$
\begin{equation}\label{e2}
{2l+M-4\choose M-4}={3M-8\choose M-4}
\end{equation}
independent conditions. Taking into account that the point $p$ and
the hyperplane $P$ are arbitrary, we get that a double cubic of
general position is regular at every point outside the
ramification divisor, provided that a violation of the condition
(R1.3) at a fixed point with a fixed hyperplane $P$ imposes on the
polynomial $g$ at least $2M$ independent conditions. It is easy to
check that for $M\geq 6$ the right-hand side of (\ref{e2}) is
higher (much higher) than $2M$. Q.E.D. for the regularity
conditions outside the ramification divisor for
$m=3$.\vspace{0.3cm}

{\bf 4.2. Proof for the case $m\geq 4$.} In that case our work
breaks into two parts: we need to check that the set
$\overline{\{q_1=q_2=0\}}\cap Q\cap P$ is either irreducible or
breaks into components in the correct way (this part of the
condition depends of the polynomial $f$ only), and, furthermore,
that the inverse image of each component with respect to $\sigma$
is irreducible (this part of the condition depends on $g$ only).
Let us start with the regularity condition on $Q$.

As we will see from the arguments below, the estimate for the set
of non-regular hypersurfaces $Q$ is the stronger, the higher gets
$m$. So it is sufficient to consider the case $m=4$. Set
$$
F_P=\overline{\{q_1=q_2=0\}}\cap P.
$$
This is a quadric of dimension $M-2$ in ${\mathbb P}^{M-1}$, more
precisely, a quadratic cone with the vertex at $p$. If $F_P$ is a
cone over a smooth quadric of dimension $M-3\geq 3$, then $F_P$ is
a factorial variety, and reducibility of the divisor
$$
(q_3+q_4)|\,_{F_P}
$$
imposes on the polynomial $q_3+q_4$ at least
$$
{M+2\choose 4}-M^2+3
$$
independent conditions, which is essentially higher than $2M$.
(This estimate is obtained for that type of reducibility, which
gives the least codimension of the non-regular set, namely, when a
section of the quadric $F_P$ by a quartic with the triple point
$p$ breaks into a hyperplane section, containing the point $p$,
and a section of $F_p$ by a cubic with the double point $p$. For
other types of reducibility the codimension is much higher.)

Therefore, we may assume that $F_P$ is a quadratic cone over a
singular quadric. Assume that $F_P$ is an irreducible quadric. In
that case $F_P$ is swept out by planes of dimension $k\geq M/2$.
Moreover, one may choose $k$-planes of general position
$L_1,L_2\subset F_P$ in such a way that their linear span
$<L_1,L_2>$ is a $(k+1)$-plane $L={\mathbb P}^{k+1}$, that is,
$L_1\cap L_2=L_{12}$ is a $(k-1)$-plane, containing the point $p$
and the vertex space of the cone $F_P$, but not coinciding with
the latter (that is, $L_{12}$ is a plane of general position,
strictly containing the vertex space of the cone $F_p$). Let us
count the independent conditions on the polynomial $(q_3+q_4)$,
which are imposed by requiring that each of the polynomials
$(q_3+q_4)|\,_{L_i}$ should be reducible. Elementary but tiresome
computations give the codimension
$$
\alpha_k=\frac{(k+5)(k+3)k(k-2)}{24}+1
$$
for the restriction of $q_3+q_4$ onto $L_i\cong{\mathbb P}^k$ to
be reducible (again, we mean the type of reducibility which gives
the least codimension, namely, into a hyperplane in ${\mathbb
P}^k$, containing the point $p$, and a cubic with the double point
$p$). Its follows from here that the required codimension for both
polynomials $(q_3+q_4)|\,_{L_i}$ to be reducible, $i=1,2$, is
$2\alpha_k-\alpha_{k-1}$, and it is not hard to check that this
integer is $\geq 2M$. For $M\geq 7$ it is sufficient to estimate
the reducibility of the polynomial $(q_3+q_4)|\,_{L_1}$,
restricted onto one $k$-plane: this already gives the required
codimension.

Finally, assume that the quadric $F_P$ is reducible, $F_P=L_1\cup
L_2$, where $L_i\cong{\mathbb P}^{M-2}$ are $(M-2)$-planes.
Reducibility of the polynomial $q_3+q_4$, restricted onto one of
the planes $L_i$, gives, in the previous notations, $\alpha_{M-2}$
independent conditions on the polynomial $q_3+q_4$, and this
number is (considerably) higher than $2M$. This completes the
proof of the $Q$-part of the condition (R1.3).

Let us check the $W$-part of the condition (R1.3). Let $Y$ be an
irreducible component of the closed set
$\overline{\{q_1=q_2=0\}}\cap Q\cap P$. We have to estimate the
number of independent conditions that follow from reducibility of
the double cover $\sigma^{-1}(Y)\to Y$. There are two possible
cases: $W\cap Y\not\subset\mathop{\rm Sing}Y$ (the case of general
position) and $W\cap Y\subset\mathop{\rm Sing}Y$. Let us treat
them separately.

If the case of general position takes place, then applying Lemma
4.1, we get at least
$$
{2l+M-4\choose M-4}\geq{M\choose 4}>2M
$$
independent conditions on the polynomial $g$, since $\mathop{\rm
dim}Y=M-3$ and $l\geq 2$. If the second case takes place, then the
hypersurface $W^*_{2l}$ contains entirely at least one irreducible
component $S$ of the set of singular points $\mathop{\rm Sing}Y$.
In that case the simplest thing to do is to apply the method of
estimating the codimension from the paper [16]: let $\pi\colon S
\to{\mathbb P}^k$, $k=M-4=\mathop{\rm dim}S$, be a generic
projection, then the pull back on $S$ of every non-trivial
polynomial on ${\mathbb P}^k$ does not vanish identically.
Therefore, the required codimension is not less than the dimension
of the space of homogeneous polynomials of degree $2l$ on
${\mathbb P}^k$, that is, again ${2l+k\choose k}\geq{M\choose
4}>2M$. This completes the proof of the $W$-part of the condition
(R1.3).

Therefore, a generic double cover $V$ is regular at every point
outside the ramification divisor.\vspace{0.3cm}

{\bf 4.3. The regularity condition on the ramification divisor.}
The fact that a generic double cover $V$ satisfies the conditions
(R2.2) and (R2.3), is proven word for word in the same way as the
conditions (R1.2) and (R1.3) were proven in the case when the
point $o$ does not lie on the ramification divisor. The condition
(R2.1) is a somewhat stronger version of the standard regularity
condition for the hypersurface $Q\cap\{z_{M+1}=0\}$ at the point
$p$: it involves an additional ($M-1$)-dimensional family of
quadratic forms $\lambda^2(z_*)-\bar{w_2}$. Let us consider it in
details.

Since the point $p$ lies on the branch divisor, to prove the
regularity condition (R2.1) it is sufficient to show that
violation of that condition imposes  on the pair of polynomials
$(f,g)$ not less than $2M-1$ conditions. For the estimates on the
rank of the quadratic forms ${\bar w}_2$ and ${\bar q}_2$ this is
easy to check. To prove the first part of the condition (R2.1),
that is, the regularity of the sequence (\ref{b2}), let us use the
following result [17].

Let $u_1,\dots,u_{N+1}$ be independent variables,
$$
{\cal P}=\prod^{k+1}_{i=1}{\cal P}_{m_i}=\{(p_1,\dots,p_{k+1})\},
$$
the space of all $(k+1)$-uples of homogeneous polynomials of
degree $m_1,\dots,m_{k+1}$, respectively, $0\leq k\leq N-1$,
$\mu=\mathop{\rm min}\{m_1,\dots,m_{k+1}\}\geq 2$. The set ${\cal
P}$ is a linear space of dimension
$$
\mathop{\rm dim}{\cal P}=\sum^{k+1}_{i=1}{m_i+N\choose N}.
$$
To every $(k+1)$-uple $(p_*)=(p_1,\dots,p_{k+1})$ we correspond
the set of its zeros $Z(p_*)\subset{\mathbb P}^N$. Thus a sequence
$p_1,\dots,p_{k+1}$ is regular in ${\cal O}_{0,{\mathbb C}^N+1}$
if and only if $\mathop{\rm codim}Z(p_*)=k+1$. Let
$$
Y=\{(p_*)\in{\cal P}|\mathop{\rm codim}\nolimits_{{\mathbb
P}^N}Z(p_*)\leq k \}
$$
be the set of non-regular sets $(p_*)$. Put $I=\{1,\dots,k+1\}$.
Set
$$
\mu_j=\mathop{\rm min}_{S\subset I,\sharp S=j}\left\{\sum_{i\in S
}m_i\right\}\geq j\mu.
$$
In [17] the following fact was shown.

{\bf Proposition 4.1.} {\it The following estimate holds:}
$$
\mathop{\rm codim}\nolimits_{\cal P}Y\geq\mathop{\rm
min}_{j\in\{0,\dots,k\}}\{(\mu_{j+1}-j)(N-j)+1\}.
$$

Applying Proposition 4.1 to our case $N=M-2$, $k+1=m$,
$(m_1,\dots,m_{k+1})=(2,2,\dots,m)$, we obtain the required
estimate for the codimension of the set of pairs $(f,g)$ that do
not satisfy the condition (R2.1). It is easy to see that this
estimate is (considerably) stronger than $2M-1$. This completes
the proof of the regularity conditions (R2.1)-(R2.3)
%%%%%%%%%%%%%%%%%%%%%%%%%%%%%%%%%%%%%%%%%%%%%%%%%%%%%%%%%%%%%%%%%
%%%%%%%%%%%%%%%%%%%%%%%%%%%%%%%%%%%%%%%%%%%%%%%%%%%%%%%%%%%%%%%%%
%%%%%%%%%%%%%%%%%%%%%%%%%%%%%%%%%%%%%%%%%Section 5

\section{Application: pencils of Fano double covers}

In this section, as an application we study movable linear systems
on Fano fiber spaces $V/{\mathbb P}^1$, the fibers of which are
double hypersurfaces. Birational geometry of those varieties was
studied in [4-6], where as a main tool the traditional quadratic
technique of the method of maximal singularities [2,18] was used,
the technique that goes back to the classical paper of
V.A.Iskovskikh and Yu.I.Manin [19]. Here we prove again the
results of [4-6] by the linear method.\vspace{0.3cm}

{\bf 5.1. Formulation of the problem and the main result.} Let
$V/{\mathbb P}^1$ be a Fano fiber space, our assumptions about
which are as follows:

(i) $V$ is a smooth projective variety with the Picard group
$\mathop{\rm Pic}V={\mathbb Z}K_V\oplus{\mathbb Z}F$, where $F$ is
the class of a fiber of the projection $\pi\colon V\to{\mathbb
P}^1$,

(ii) every fiber $F_t=\pi^{-1}(t)$, $t\in{\mathbb P}^1$, is a Fano
double hypersurface, $F_t\in{\cal F}$, and moreover, $F_t$ is
regular at every smooth point (in particular, every smooth fiber
$F_t\in{\cal F}_{\rm reg}$),

(iii) if $o\in F_t$ is a singular point, then $o$ is an isolated
quadratic singularity (in particular, there are finitely many such
points), satisfying the corresponding regularity condition (R1.4)
or (R2.4), depending on whether the point $o$ lies outside the
ramification divisor or on that divisor; the conditions (R1.4) and
(R2.4) are formulated below.

{\bf Theorem 3.} (i) {\it Assume that $\Sigma\subset|-nK_V+lF|$ is
a movable linear system on $V$ and $l\in{\mathbb Z}_+$. Then the
following equality holds:
$$
c_{\rm virt}(\Sigma)=c(\Sigma)=n.
$$
{\rm (ii)} Assume in addition that the Fano fiber space
$V/{\mathbb P}^1$ satisfies the $K$-condition:
$$
-K_V\not\in\mathop{\rm Int}A^1_{\rm mov}V,
$$
where $A^1_{\rm mov}V\subset A^1_{\mathbb R} V\cong {\mathbb R}^2$
is the closed cone generated by the classes of movable divisors on
$V$. Then the fiber space $V/{\mathbb P}^1$ is birationally
superrigid, in particular, the projection $\pi\colon V\to {\mathbb
P}^1$ is the only non-trivial structure of a rationally connected
fiber space on $V$.}

{\bf Proof.} It is sufficient to check the first claim: the second
one follows from it in a straightforward way. Assume that $c_{\rm
virt}(\Sigma)<c(\Sigma)$, then the linear system $\Sigma$ has a
maximal singularity $E\subset V^+$, where $\varphi\colon V^+\to V$
is a birational morphism, $E$ is an exceptional divisor. Its
centre $B=\varphi(E)\subset V$ is an irreducible subvariety of
codimension 2 or higher. Let $F=F_t$ be some fiber, intersecting
$B$. If $B\cap F$ is not a singular point of the fiber $F$, then,
restricting the system $\Sigma$ onto the fiber $F$, we get a
contradiction with Theorem 1: for a generic divisor $D\in\Sigma$
the pair $(F,\frac{1}{n}D|\,_F)$ is not canonical. Therefore, one
may assume that $B=o\in F$ is a singular point of a fiber. As we
will show below in Sec. 5.2-5.3, this possibility does not
realize, either, if the corresponding regularity condition (R1.4)
or (R2.4) holds. Therefore, the linear system $\Sigma$ cannot have
maximal singularities, which proves the theorem.

{\bf Remark 5.1.} As the proof given above shows, the linear
method does not make use of the condition of twistedness of the
fiber space $V/{\mathbb P}^1$ over the base. Essentially, all the
proof is restricting to some fiber and applying the property of
divisorial (log) canonicity of every fiber. All work is
accumulated in the proof of that property, which requires stronger
regularity conditions than birational rigidity. The regularity
conditions, on which the quadratic method is based, are weaker and
for that reason the results obtained by the quadratic method, are
more precise. However, the quadratic constructions require
considerably more work.\vspace{0.3cm}

{\bf 5.2. Singular points of fibers.} Let  $o\in F=F_t$ be a
singularity of the fiber $\pi^{-1}(t)$. We say that the point $o$
is a singularity {\it of the first type}, if it lies outside the
ramification divisor of the morphism $\sigma\colon F\to
G\subset{\mathbb P}$, that is, $p=\sigma(o)\not\in W$. In that
case we assume that the hypersurface $G$ in suitable affine
coordinates with the origin at the point $p$ is given by the
equation
$$
f=q_2+\dots+q_m=0,
$$
where $q_2(z_1,\dots,z_{M+1})$ is a non-degenerate quadratic form.
The branch divisor of the morphism $\sigma$ is cut out on $G$ by
the hypersurface
$$
g=1+w_1+\dots+w_{2l}=0,
$$
which does not pass through $p$. We say that the point $o$ is a
singularity {\it of the second type}, if it lies on the
ramification divisor of the morphism $\sigma$, that is,
$p=\sigma(o)\in W$. In that case $G$ is given by the equation
$$
f=q_1+\dots+q_m=0,
$$
where $q_1\equiv z_{M+1}$, and the equation of the branch divisor
is of the form
$$
g=\alpha z_{M+1}+w_2+\dots+w_{2l}=0,
$$
where $\alpha\in{\mathbb C}$ is an arbitrary constant.

Let us consider first the singularities of the first type. We say
that the variety $F$ is {\it regular} at the point $o$ of the
first type, if (in the notations above) the following condition
holds:

(R1.4) for any non-zero linear form $\lambda(z_1,\dots,z_{M+1})$
for $m\leq 2l$ the sequence
$$
\lambda,q_2,\dots,q_m,h_{l+1},\dots,h_{2l-1}
$$
is regular in the local ring ${\cal O}_{0,{\mathbb C}^{M+1}}$, and
for $m>2l$ the sequence
$$
\lambda,q_2,\dots,q_{m-1},h_{l+1},\dots,h_{2l}
$$
is regular in the ring ${\cal O}_{0,{\mathbb C}^{M+1}}$, where
$h_j(z_*)$ have the same meaning as in the regularity condition
(R1.1) and, moreover, the closed set
$\sigma^{-1}(\{\lambda=q_2=0\}\cap G)$ is irreducible.

Since there are finitely many singular points, the regularity
condition does not require a special proof (it is sufficient to
require that the system of equations
$$
q_2=\dots=q_m=h_{l+1}=\dots=h_{2l-1}=0
$$
or, respectively, the one with $q_m$ replaced by $h_{2l}$,
determines a closed set, none of the components of which is
contained in a hyperplane; the second part of the regularity
condition holds at a point of general position in an obvious way).

Let $D\in|-nK_F|$ be a prime divisor, $\beta\colon F^+\to F$ the
blow up of the point $o$, $\beta_G\colon G^+\to G$ the blow up of
the point $p$ on $G$, $E_F\subset F^+$ and $E_G\subset G^+$ the
exceptional divisors (non-singular quadrics), where $E_F$
identifies naturally with $E_G$. Let $B\subset E_F$ be a
hyperplane section. For the strict transform $D^+\subset F^+$ we
have $D^+\sim-nK_F-\nu E_F$, where $\nu\in{\mathbb Z}_+$ is some
integer.

{\bf Proposition 5.1.} {\it The following inequality holds:}
\begin{equation}\label{f1}
\nu+\mathop{\rm mult}\nolimits_BD^+\leq 2n.
\end{equation}

{\bf Proof.} Assume the converse. Let $R\subset{\mathbb P}$ be the
unique hyperplane, $R\ni p$, cutting out $B$ on $E_G\cong E_F$,
say, $R=\{\lambda(z_*)=0\}$. The divisor $T=\sigma^{-1}(R\cap G)$
satisfies the inequality (\ref{f1}) with $n=1$ and is irreducible,
so that $D\neq T$ and one may form the scheme-theoretic
intersection $Y=(D\circ T)$, an effective cycle of codimension 2,
satisfying the inequality
$$
\frac{\mathop{\rm mult}_o}{\mathop{\rm deg}}(D\circ
T)>\frac{4}{\mathop{\rm deg}F}=\frac{2}{m}.
$$
By linearity, we may assume that the cycle $Y$ is an irreducible
subvariety of codimension 2 on $F$, or a prime divisor on $T$. Set
$D_2$ to be the divisor on $F$, cut out by the hypersurface
$q_2=0$. By the regularity condition, $(D_2\circ T)$ is a prime
divisor on $T$, and moreover,
$$
\frac{\mathop{\rm mult}_o}{\mathop{\rm deg}}(D_2\circ
T)=\frac{3}{2m}.
$$
Therefore, $Y\neq(D_2\circ T)$, so that $Y\not\subset D_2$ and we
may form the effective cycle $(Y\circ D_2)$ that has an
irreducible component $Y^{\sharp}$, satisfying the inequality
$$
\frac{\mathop{\rm mult}_o}{\mathop{\rm
deg}}Y^{\sharp}>\frac{3}{m}.
$$
The subvariety $Y^{\sharp}$ is of codimension two on $T$. The
condition (R1.4) can be understood as the regularity condition for
the variety $T$. Applying the technique of hypertangent divisors
on $T$ to the subvariety $Y^{\sharp}$ in the standard way, we
obtain a contradiction (intersecting one by one with hypertangent
divisors, we construct a curve $C\subset T$, satisfying the
inequality $(\mathop{\rm mult}_o/\mathop{\rm deg})C>1$, which is
impossible). Q.E.D. for Proposition 5.1. The case of a singularity
of the first type is excluded.\vspace{0.3cm}

{\bf 5.3. Singularities of the second type.} We say that the
variety $F$ is {\it regular} at a point $o$ of the second type, if
the following condition holds.

(R2.4) For any linear form $\lambda(z_1,\dots,z_M)$ the sequence
$$
\lambda^2(z_*)-{\bar w}_2,{\bar q}_2,\dots,{\bar q}_m
$$
is regular in ${\cal O}_{0,{\mathbb C}^{M+1}}$, where ${\bar
w}_i,{\bar q}_j$ stand for the restrictions of the polynomials
$w_i,q_j$ onto the hyperplane $z_{M+1}=0$, whereas the quadratic
forms ${\bar w}_2$ and ${\bar q}_2$ have the full rank, the closed
set $\{{\bar q}_2={\bar q}_3=0\}$ in ${\mathbb C}^M$ is
irreducible and its linear span is the whole space ${\mathbb
C}^M$, whereas the closed set
$$
\sigma^{-1}(\overline{\{\lambda(z_*)=q_1=q_2=0\}}\cap G)
$$
is irreducible and reduced.

Note again, that since there are finitely many singular points,
the condition (R2.4) does not need a special proof.

Let $D\in|-nK_F|$ be a prime divisor, $\beta\colon F^+\to F$ and
$\beta_G\colon G^+\to G$ the blow ups of the points $o$ and
$p=\sigma(o)$, respectively, with the exceptional divisors
$E_F\subset F^+$ and $E_G\subset G^+$. The morphism $\sigma$
extends to the double cover $\sigma^+\colon F^+\to G^+$, where
$\sigma_E\colon E_F\to E_G$ realizes the quadric $E_F$ as the
double cover of $E_G\cong{\mathbb P}^{M-1}_{(z_1:\dots:z_m)}$
branched over $\{{\bar w}_2=0\}$. Let $D^+\subset F^+$ be the
strict transform of the divisor $D$ and define the integer
$\nu\in{\mathbb Z}_+$ by the equivalence $D^+\sim-nK_F-\nu E_F$.
Let $B\subset E_F$ be a hyperplane section of the smooth quadric
$E_F$. Repeating the arguments of \S 2 word for word (considering
the cases when $B$ is pulled back and not pulled back from $G$
separately), we obtain the claim of Proposition 5.1. This excludes
the case of a singularity of the second type. Proof of Theorem 3
is complete.

%%%%%%%%%%%%%%%%%%%%%%%%%%%%%%%%%%%%%%%%%%%%%%%%%%%%%%%%%%
%\newpage

{\small

\section*{References}

\noindent 1. A.V.Pukhlikov, Birational geometry of Fano direct
products, Izvestiya: Mathematics, V. 69 (2005), no. 6,
1225-1255.\vspace{0.3cm}

\noindent 2. Pukhlikov A.V., Birationally rigid varieties. I. Fano
varieties. Russian Math. Surveys. V. 62 (2007), No. 5.
\vspace{0.3cm}

\noindent 3. Pukhlikov A.V., Birationally rigid Fano double
hypersurfaces, Sbornik: Mathematics {\bf 191} (2000), No. 6,
883-908. \vspace{0.3cm}

\noindent 4. A.V.Pukhlikov, Birationally rigid varieties with a
pencil of Fano double covers. III,  Sbornik: Mathematics. V. 197
(2006), no. 3, 335-368.\vspace{0.3cm}

\noindent 5. A.V.Pukhlikov, Birationally rigid varieties with a
pencil of Fano double covers. I,  Sbornik: Mathematics. V. 195
(2004), no. 7, 1039-1071.\vspace{0.3cm}

\noindent 6. A.V.Pukhlikov, Birationally rigid varieties with a
pencil of Fano double covers. II,  Sbornik: Mathematics. V. 195
(2004), no. 11, 1665-1702.\vspace{0.3cm}

\noindent 7. Pukhlikov A.V., Birational geometry of algebraic
varieties with a pencil of Fano complete intersections,
Manuscripta Mathematica V.121 (2006), 491-526. \vspace{0.3cm}

\noindent 8. Corti A., Singularities of linear systems and 3-fold
birational geometry, in ``Explicit Birational Geometry of
Threefolds'', London Mathematical Society Lecture Note Series {\bf
281} (2000), Cambridge University Press, 259-312. \vspace{0.3cm}

\noindent 9. Cheltsov I.A. Birationally rigid Fano varieties.
Russian Math. Surveys. V. 60 (2005), no. 5, 875-965.\vspace{0.3cm}

\noindent 10. Cheltsov I.A. and Grinenko M.M., Birational rigidity
is not an open property, preprint, 2006, arXiv:math.AG/0612159.
\vspace{0.3cm}

\noindent 11. Tian G., On K\" ahler-Einstein metrics on certain
K\" ahler manifolds with $c_1(M)>0$. Invent. Math. {\bf 89}
(1987), 225-246. \vspace{0.3cm}

\noindent 12. Demailly J.-P. and Koll\' ar J., Semi-continuity of
complex singularity exponents and K\" ahler-Einstein metrics on
Fano orbifolds, Ann. Sci. de l'{\' E}cole Norm. Sup. {\bf 34}
(2001), 525-556. \vspace{0.3cm}

\noindent 13. Nadel A., Multiplier ideal sheaves and K\"
ahler-Einstein metrics of positive scalar curvature. Ann. Math.
{\bf 132} (1990), 549-596. \vspace{0.3cm}

\noindent 14. Pukhlikov A.V., Birationally rigid iterated Fano
double covers. Izvestiya: Mathematics. {\bf 67} (2003), no. 3,
555-596. \vspace{0.3cm}

\noindent 15. Fulton W., Intersection Theory, Springer-Verlag,
1984. \vspace{0.3cm}

\noindent 16. Pukhlikov A.V., Birational automorphisms of Fano
hypersurfaces, Invent. Math. {\bf 134} (1998), no. 2, 401-426.
\vspace{0.3cm}

\noindent 17. Pukhlikov A.V., Birationally rigid Fano complete
intersections, Crelle J. f\" ur die reine und angew. Math. {\bf
541} (2001), 55-79. \vspace{0.3cm}

\noindent 18. Pukhlikov A.V., Birational automorphisms of
three-dimensional algebraic varieties with a pencil of del Pezzo
surfaces, Izvestiya: Mathematics {\bf 62}:1 (1998),
115-155.\vspace{0.3cm}

\noindent 19. Iskovskikh V.A. and Manin Yu.I., Three-dimensional
quartics and counterexamples to the L\" uroth problem, Math. USSR
Sb. {\bf 86} (1971), no. 1, 140-166.%\vspace{0.3cm}

}

\begin{flushleft}
e-mail: {\it pukh@liv.ac.uk, pukh@mi.ras.ru}
\end{flushleft}

\end{document}